\overfullrule=0pt
\centerline {\bf MISCELLANOUS APPLICATIONS OF CERTAIN MINIMAX THEOREMS. I}\par
\bigskip
\bigskip
\centerline {B. RICCERI}\par
\bigskip
\bigskip
\centerline {Department of Mathematics, University of Catania}\par
\centerline {Viale A. Doria 6, 95125 Catania, Italy}
\bigskip
\bigskip
{\bf ABSTRACT.} Here is one of the results of this paper (with the convention ${{1}\over {0}}=+\infty$): Let $X$ be a real Hilbert space and 
let $J:X\to {\bf R}$ be a $C^1$ functional, with compact derivative, such that 
$$\alpha^*:=\max\left \{0,\limsup_{\|x\|\to +\infty}{{J(x)}\over {\|x\|^2}}\right \}<\beta^*:=\sup_{x\in X\setminus \{0\}}{{J(x)}\over 
{\|x\|^2}}<+\infty\ .$$
Then, for every $\lambda\in \left ]{{1}\over {2\beta^*}}, {{1}\over {2\alpha^*}}\right [$ and
for every convex set $C\subseteq X$ dense in $X$, there exists $\tilde y\in C$ such that the equation
$$x=\lambda J'(x)+\tilde y$$
has at least three solutions, two of which are global minima of the functional
$x\to {{1}\over {2}}\|x\|^2-\lambda J(x)-\langle x,\tilde y\rangle$\ .\par
\bigskip
\bigskip
{\bf AMS (MOS) Subject Classification.} 49K35; 90C47; 47G40; 35A15; 35J20
\bigskip
\bigskip
\bigskip
\bigskip
\centerline {1. {\bf RESULTS}}
\bigskip
\bigskip
The present one is the first of a series of notes (with the same title) on new consequences and
applications of certain minimax theorems that we have established in the past years ([2]-[10]).
\smallskip
In [8], in particular, we obtained the following result:
\medskip
{\bf Theorem A.} ([8], Theorem 3.4) {\it Let $X$ be a non-empty set, $Y$  a real inner product space and
$I:X\to {\bf R}$, $\Phi:X\to Y$ two given functions.\par
Then, for each $\mu>0$, at least one of the following assertions holds:\par
\noindent
$(a)$\hskip 5pt for each filtering cover ${\cal N}$ of $X$, there exists
$A\in {\cal N}$ such that
$$\sup_{y\in Y}\inf_{x\in A}(I(x)+\mu(2\langle\Phi(x),y\rangle-
\|y\|^2))<
\inf_{x\in A}\sup_{y\in \Phi(A)}(I(x)+\mu(2\langle\Phi(x),y\rangle-\|y\|^2))\ ;$$
\noindent
$(b)$\hskip 5pt for each global minimum $u$ of $x\to I(x)+\mu\|\Phi(x)\|^2$,
one has
$$I(u)\leq I(x)+2\mu(\langle\Phi(x),\Phi(u)\rangle-\|\Phi(u)\|^2)$$
for all $x\in X$.}\par
\medskip
A cover ${\cal N}$ of the set $X$ is said to be filtering if for each $A_1, A_2\in {\cal N}$, there
exists $A_3\in {\cal N}$ such that $A_1\cup A_2\subseteq A_3$.\par
\smallskip
In the present short note, we want to highlight two applications of the following consequences of Theorem A:\par
\medskip
{\bf Theorem 1.1.}  {\it Let $X$ be a non-empty set, $x_0\in X$, $Y$ a real inner product space,
$I:X\to {\bf R}$, $\Phi:X\to Y$, with $I(x_0)=0$, $\Phi(x_0)=0$, and $\mu>0$. Assume that
$$\inf_{x\in X}I(x)<0\leq\inf_{x\in X}(I(x)+\mu\|\Phi(x)\|^2)\ .$$
Then,  for each filtering cover ${\cal N}$ of $X$, there exists
$A\in {\cal N}$ such that
$$\sup_{y\in Y}\inf_{x\in A}(I(x)+\mu(2\langle\Phi(x),y\rangle-
\|y\|^2))<
\inf_{x\in A}\sup_{y\in \Phi(A)}(I(x)+\mu(2\langle\Phi(x),y\rangle-\|y\|^2))\ .$$}\par
\smallskip
{\it Proof}. The assumptions imply that $x_0$ is a global minimum of $x\to I(x)+\mu\|\Phi(x)\|^2$. But, at the same time,
since $\inf_XI<0$, $x_0$ is not a global minimum of $I$. Hence, $(b)$ of Theorem A does not hold and so $(a)$ holds. \hfill
$\bigtriangleup$
\medskip
{\bf Theorem 1.2.}  {\it Let $X$ be a non-empty symmetric set in a real vector space, $Y$ a real inner product space,
$I:X\to {\bf R}$ an even function and $\Phi:X\to Y$ an odd function.\par
Then, for each $\mu>0$, at least one of the following assertions holds:\par
\noindent
$(a_1)$\hskip 5pt for each filtering cover ${\cal N}$ of $X$, there exists
$A\in {\cal N}$ such that
$$\sup_{y\in Y}\inf_{x\in A}(I(x)+\mu(2\langle\Phi(x),y\rangle-
\|y\|^2))<
\inf_{x\in A}\sup_{y\in \Phi(A)}(I(x)+\mu(2\langle\Phi(x),y\rangle-\|y\|^2))\ ;$$
\noindent
$(b_1)$\hskip 5pt if $u\in X$ is a global minimum of $x\to I(x)+\mu\|\Phi(x)\|^2$, then $\Phi(u)=0$ and $u$ is a global minimum of $I$.}\par
\smallskip
{\it Proof}. Assume that $(a_1)$ does not hold. Let $u\in X$ be any global minimum of  $x\to I(x)+\mu\|\Phi(x)\|^2$. Then, by Theorem A,
we have
$$I(u)\leq I(x)+2\mu(\langle\Phi(x),\Phi(u)\rangle-\|\Phi(u)\|^2)\eqno{(1)}$$
for all $x\in X$. So, in particular, we have
$$I(u)\leq I(-u)+2\mu(\langle \Phi(-u),\Phi(u)\rangle-\|\Phi(u)\|^2)= I(u)-4\mu\|\Phi(u)\|^2$$
which yields
$$\|\Phi(u)\|^2=0\ .$$
So, by $(1)$, we get
$$I(u)\leq I(x)$$
for all $x\in X$.\hfill $\bigtriangleup$\par
\medskip
The application of Theorem 1.1 we wish to present is provided by the following result jointly with Theorem 1.4:\par
\medskip
{\bf Theorem 1.3.}  {\it Let $X$ be a real inner product space and let $\tau$ be a topology on $X$.
Moreover, let $J:X\to {\bf R}$ be a functional such that
$$J(0)=0<\sup_X J$$
and
$$\beta^*:=\sup_{x\in X\setminus \{0\}}{{J(x)}\over {\|x\|^2}}<+\infty\ .\eqno{(2)}$$
Finally, let $\lambda>{{1}\over {\beta^*}}$ and let ${\cal N}$ be a filtering cover of $X$ such that, for
each $A\in {\cal N}$ and each $y\in X$, the restriction to $A$ of the functional $x\to \|x\|^2-\lambda J(x)+\langle x,y\rangle$ is
$\tau$-lower semicontinuous and inf-$\tau$-compact.\par
Then, there exists $\tilde A\in {\cal N}$ with the following property: for every convex set $C\subseteq X$ whose closure (in the strong 
topology) contains $\tilde A$, there exists $\tilde y\in C$ such that the restriction to $\tilde A$ of the functional 
$x\to \|x\|^2-\lambda J(x)+\left \langle x,2(\beta^*\lambda-1)\tilde y\right \rangle$ has
at least two global minima.}\par
\smallskip
{\it Proof}. In view of $(2)$, we have 
$$\inf_{x\in X}(\|x\|^2-\lambda J(x))<0 \ ,$$
as well as
$$\inf_{x\in X}(\|x\|^2-\lambda J(x)+(\beta^*\lambda-1)\|x\|^2)\geq 0\ .$$
So, we can apply Theorem 1.1 taking $Y=X$,
$$\mu=\beta^*\lambda-1\ ,$$
$$I(x)=\|x\|^2-\lambda J(x)$$
and 
$$\Phi(x)=x\ .$$
Therefore, there exists $\tilde A\in {\cal N}$ such that
$$\sup_{y\in Y}\inf_{x\in \tilde A}(\|x\|^2-\lambda J(x)+(\beta^*\lambda-1)(2\langle x,y\rangle-
\|y\|^2))<
\inf_{x\in \tilde A}\sup_{y\in \tilde A}(\|x\|^2-\lambda J(x)+(\beta^*\lambda-1)(2\langle x,y\rangle-\|y\|^2))\ .\eqno{(3)}$$\par
Now, consider the function $f:X\times X\to {\bf R}$ defined by
$$f(x,y)=\|x\|^2-\lambda J(x)+(\beta^*\lambda-1)(2\langle x,y\rangle-\|y\|^2)$$
for all $(x,y)\in X\times X$. 
 Since $f(x,\cdot)$ is continuous and $\tilde A\subseteq \overline {C}$,
 we have
$$\sup_{y\in \tilde A}f(x,y)=\sup_{v\in \overline {\tilde A}}f(x,y)\leq \sup_{v\in \overline {C}}f(x,y)
= \sup_{y\in C}f(x,y)$$
for all $x\in X$, and hence, taking $(3)$ into account, it follows that 
$$\sup_{y\in C}\inf_{x\in \tilde A}f(x,y)<
\inf_{x\in \tilde A}\sup_{y\in \tilde A}f(x,y)\leq \inf_{x\in \tilde A}\sup_{y\in C}f(x,y)\ .\eqno{(4)}$$
Now, in view of $(4)$, taking into account that $f_{|\tilde A\times C}$ is $\tau$-lower semicontinuous and inf-$\tau$-compact in $\tilde A$, and continuous
and concave in $C$, we can apply Theorem 3.2 of [8] to $f_{|\tilde A\times C}$. Consequently, there exists $\tilde y\in C$ such that
$f_{|\tilde A}(\cdot,\tilde y)$ has at least two global minima, and the proof is complete. \hfill $\bigtriangleup$\par
\medskip
From Theorem 1.3, in turn, we obtain the following result (with the convention ${{1}\over {0}}=+\infty$):\par
\medskip
{\bf Theorem 1.4.} {\it Let $X$ be a real Hilbert space and let $J:X\to {\bf R}$ be a $C^1$ functional, with compact derivative,
such that 
$$\alpha^*:=\max\left \{0,\limsup_{\|x\|\to +\infty}{{J(x)}\over {\|x\|^2}}\right \}<\beta^*:=\sup_{x\in X\setminus \{0\}}{{J(x)}\over 
{\|x\|^2}}<+\infty\ .$$
Then, for every $\lambda\in \left ]{{1}\over {2\beta^*}}, {{1}\over {2\alpha^*}}\right [$ and
for every convex set $C\subseteq X$ dense in $X$, there exists $\tilde y\in C$ such that the equation
$$x=\lambda J'(x)+\tilde y$$
has at least three solutions, two of which are global minima of the functional
$x\to {{1}\over {2}}\|x\|^2-\lambda J(x)-\langle x,\tilde y\rangle$\ .}\par
\smallskip
{\it Proof}. Fix $\lambda\in \left ]{{1}\over {2\beta^*}}, {{1}\over {2\alpha^*}}\right [$ and
a convex set $C\subseteq X$ dense in $X$. For each $y\in X$, we have
$$\liminf_{\|x\|\to +\infty}\left ( 1-2\lambda{{J(x)}\over {\|x\|^2}}-{{\langle x,y\rangle}\over {\|x\|^2}}\right )=
1-2\lambda\limsup_{\|x\|\to +\infty}{{J(x)}\over {\|x\|^2}}>0\ .$$
So, from the identity
$$\|x\|^2-2\lambda J(x)-\langle x,y\rangle=\|x\|^2\left ( 1-2\lambda{{J(x)}\over {\|x\|^2}}-{{\langle x,y\rangle}\over {\|x\|^2}}\right )$$
it follows that
$$\lim_{\|x\|\to +\infty} (\|x\|^2-2\lambda J(x)-\langle x,y\rangle)=+\infty\ .\eqno{(5)}$$
Since $J'$ is compact, $J$ is sequentially weakly continuous ([12], Corollary 41.9).  Then, in view of $(5)$ and of the Eberlein-Smulyan theorem, for each
$y\in X$, the functional $x\to \|x\|^2-2\lambda J(x)+\langle x,y\rangle$ is inf-weakly compact in $X$. So, we can apply Theorem 1.3 taking the weak topology as $\tau$
and ${\cal N}=\{X\}$. Consequently, since the set ${{1}\over {1-2\beta^*\lambda}}C$ is convex and dense in $X$, there exists $\hat y\in {{1}\over {1-2\beta^*\lambda}}C$ such that the functional  $x\to \|x\|^2-2\lambda J(x)+\left \langle x,2(2\beta^*\lambda-1)\hat y\right \rangle$ has
at least two global minima in $X$ which are two of its critical points. Since the same functional satisfies the Palais-Smale condition ([12], Example 38.25), it has a third critical 
point in view of Corollary 1 of [1]. Clearly, the conclusion follows taking $\tilde y=(1-2\beta^*\lambda)\hat y$\ .\hfill $\bigtriangleup$\par
\medskip
We now give a specific application of Theorem 1.4.\par
\smallskip
Let $\Omega\subset {\bf R}^n$ be a bounded domain with smooth
boundary. On the Sobolev space $H^1_0(\Omega)$, we
consider the scalar product
$$\langle u,v\rangle=\int_{\Omega}\nabla u(x)\nabla v(x)dx$$
with the induced norm
$$\|u\|=\left ( \int_{\Omega}|\nabla u(x)|^2 dx\right ) ^{1\over 2}\ .$$
We denote by $H^{-1}(\Omega)$ the dual of $H^1_0(\Omega)$.\par
\smallskip
If $n\geq 2$, we denote by ${\cal A}$ the class of all
Carath\'eodory functions $f:\Omega\times {\bf R}\to {\bf R}$ such that
$$\sup_{(x,\xi)\in \Omega\times {\bf R}}{{|f(x,\xi)|}\over
{1+|\xi|^q}}<+\infty\ ,$$
where  $0<q< {{n+2}\over {n-2}}$ if $n>2$ and $0<q<+\infty$ if
$n=2$. While, when $n=1$, we denote by ${\cal A}$  the class
of all Carath\'eodory functions $f:\Omega\times {\bf R}\to {\bf R}$ such
that, for each $r>0$, the function $x\to \sup_{|\xi|\leq r}|f(x,\xi)|$ belongs
to $L^{1}(\Omega)$.\par
\smallskip
Given $f\in {\cal A}$ and $\varphi\in H^{-1}(\Omega)$,  consider the following Dirichlet problem
$$\cases {-\Delta u= f(x,u) +\varphi
 & in
$\Omega$\cr & \cr u=0 & on
$\partial \Omega$\ .\cr}\eqno{(P_{f,\varphi})} $$
 Let us recall
that a weak solution
of $(P_{f,\varphi})$ is any $u\in H^1_0(\Omega)$ such that
 $$\int_{\Omega}\nabla u(x)\nabla v(x)dx
-\int_{\Omega}f(x,u(x))v(x)dx-\varphi(v)=0$$
for all $v\in H^1_0(\Omega)$.\par
\smallskip
Let $\Psi, J_{f}:H^1_0(\Omega)\to {\bf R}$ be the functionals defined by
$$\Psi(u)={{1}\over {2}}\|u\|^2$$
$$J_{f}(u)=\int_{\Omega}F(x,u(x))dx\ ,$$
where
$$F(x,\xi)=\int_{0}^{\xi}f(x,t)dt\ .$$
Notice that $J_f$ is a $C^1$ functional whose derivative is given by
$$J'_{f}(u)(v)=\int_{\Omega}f(x,u(x))v(x)dx$$
for all $u,vH^1_0(\Omega)$. Consequently, the weak solutions 
of problem $(P_{f,\varphi})$ are exactly the critical points in
$H^1_0(\Omega)$ of the functional $\Psi-J_{f}-\varphi$. Moreover,
$J'_{f}$ is compact. \par
\smallskip
Furthermore, if $f\in {\cal A}$ and the set 
$$\left\{x\in \Omega : \sup_{\xi\in {\bf R}}F(x,\xi)>0\right\}$$
has a positive measure, we have $\sup_{H^1_0(\Omega)}J_f>0$ (see the proof of
Theorem 2 of [9]).\par
\smallskip
We now state the following\par
\medskip
{\bf Theorem 1.5.} {\it Let $f\in {\cal A}$ be such that the set
$$\left\{x\in \Omega : \sup_{\xi\in {\bf R}}F(x,\xi)>0\right\}$$
has a positive measure and
$$\limsup_{|\xi|\to +\infty}{{\sup_{x\in\Omega}F(x,\xi)}\over {\xi^2}}\leq 0\ .\eqno{(6)}$$
Then, for every $\lambda>0$ large enough and for every convex set $C\subset H^{-1}(\Omega)$ dense
in $H^{-1}(\Omega)$, there exists $\varphi\in C$ such that
the problem
$$\cases {-\Delta u= \lambda f(x,u) +\varphi
 & in
$\Omega$\cr & \cr u=0 & on
$\partial \Omega$\cr} $$
has at least three weak solutions, two of which are global minima in $H^1_0(\Omega)$ of the functional $\Psi-\lambda J_f-\varphi$.}\par
\smallskip
{\it Proof}. Condition $(6)$ clearly implies that
$$\limsup_{\|u\|\to +\infty}{{J_f(u)}\over {\|u\|^2}}\leq 0\ .$$
Then, in view of the above-recalled preliminaries, we can apply Theorem 1.4 taking $X=H^{1}_0(\Omega)$,  $J=J_f$, and the conclusion
directly follows.\hfill $\bigtriangleup$
\par
\medskip
The last result is an application of Theorem 1.2.\par
\smallskip
Let $(T,{\cal F},\mu)$ be a measure space, 
 $E$  a real Banach space and $p\geq 1$.
\par
\smallskip
As usual, $L^{p}(T,E)$  denotes the space of all (equivalence
classes of) strongly $\mu$-measurable functions $u : T\rightarrow E$ 
such that
$\int_{T}\parallel u(t)\parallel^{p} d\mu<+\infty$, equipped with
the norm $$\parallel u\parallel_{L^{p}(T,E)}=
\left ( \int_{T}\parallel u(t)\parallel^{p}d\mu\right ) ^{1\over p}\ .$$
\smallskip
A set $D\subseteq L^{p}(T,E)$ is said to be decomposable if, for
every $u,v\in D$ and every $A\in {\cal F}$, the function
 $$t\to
\chi_{A}(t)u(t)+(1-\chi_{A}(t))v(t)$$ belongs to $D$, where $\chi_{A}$
denotes the characteristic function of $A$.\par
\medskip
{\bf Theorem 1.6.}  {\it  Let $(T,{\cal F},\mu)$ be a non-atomic
measure space, with $0<\mu(T)<+\infty$, and $E$ a real Banach space. Let $f:E\to {\bf R}$ be
a lower semicontinuous and even function which has no global minima. Let $g:T\times E\to {\bf R}$ be a Carath\'eodory
function such that $g(x,\cdot)$ is odd for all $t\in T$. Moreover, assume that
$$\max\left \{ \sup_{x\in E}{{|f(x)|}\over {1+\|x\|^p}},
\sup_{(t,x)\in T\times E}{{|g(t,x)|}\over {1+\|x\|^p}}
\right \} <+\infty$$
for some $p\geq 1$.
For each $u\in L^p(T,E)$, set
$$J(u)=\int_Tf(u(t))d\mu+\left ( \int_Tg(t,u(t))d\mu\right )^2\ .$$
Then, the restriction of
the functional $J$
to any symmetric and decomposable subset of $L^p(T,E)$ containing the constant
functions has no global
minima.}\par
\smallskip
{\it Proof}.  Let $X\subseteq L^p(T,E)$ be a symmetric and decomposable set containing
the constant functions. Arguing by contradiction, assume that
$\hat u\in X$ is a global minimum of $J_{|X}$. For each $h\in L^p(T)$,
set
$$X_h=\{u\in X : \|u(t)\|\leq h(t)\hskip 3pt \hbox {\rm a.e.\hskip 3pt
in\hskip 3pt T}\}\ .$$
Clearly, $X_h$ is decomposable and
the family $\{X_h\}_{h\in L^p(T)}$ is a filtering covering of
$X$. Let $I, \Phi:L^p(T,E)\to {\bf R}$ be the functionals defined by
$$I(u)=\int_Tf(u(t))d\mu$$
and
$$\Phi(u)=\int_Tg(t,u(t))d\mu$$
for all $u\in L^p(T,E)$.
Of course, $I$ is lower semicontinuous and even, while $\Phi$ is continuous and odd.
Fix $h\in L^p(T)$. Notice that $X_h$ is connected ([11]). Therefore, $\Phi(X_h)$
is an interval, and so 
$$\Lambda_h:=\overline {\Phi(X_h)}$$
is a compact interval. Now,
we can apply Theorem 1.D of [6]. Thanks to it, we then have
$$\sup_{\lambda\in \Lambda_h}\inf_{u\in X_h}(I(u)+2\lambda\Phi(u)-\lambda^2)=
\inf_{u\in X_h}\sup_{\lambda\in \Lambda_h}(I(u)+2\lambda\Phi(u)-\lambda^2)\ .$$
Hence, (with $\mu=1$) assertion $(a_1)$ of Theorem 1.2 does not hold, and so
assertion $(b_1)$ must hold. 
Therefore,  $\hat u$ is a global minimum of $I_{|X}$ and $\Phi(\hat u)=0$. 
Since $X$ contains the constants, we have
$$I(\hat u)\leq \mu(T)\inf_E f\ .$$
This clearly implies that
$$\int_T(f(\hat u(t))-\inf_Ef)d\mu=0\ ,$$
and so $\hat u(t)$ would be a global minimum of $f$ for $\mu$-almost $t\in T$,
against the assumptions.\hfill $\bigtriangleup$
\bigskip
\bigskip

{\bf Acknowledgement.} The author has been supported by the Gruppo Nazionale per l'Analisi Matematica, la Probabilit\`a 
e le loro Applicazioni (GNAMPA) of the Istituto Nazionale di Alta Matematica (INdAM).

\vfill\eject
\centerline {\bf REFERENCES}\par
\bigskip
\bigskip
\noindent
[1]\hskip 5pt P. Pucci and J. Serrin, A mountain pass theorem,
{\it J. Differential Equations}, 60:142-149, 1985.\par
\smallskip
\noindent
[2]\hskip 5pt B. Ricceri, Some topological mini-max theorems via
an alternative principle for multifunctions, {\it Arch. Math.} (Basel),
 60:367-377, 1993.\par
\smallskip
\noindent
[3]\hskip 5pt B. Ricceri, On a topological minimax theorem and
its applications, in ``Minimax theory and applications'', B. Ricceri
and S. Simons eds., Kluwer Academic Publishers, 191-216, 1998.\par
\smallskip
\noindent
[4]\hskip 5pt B. Ricceri, A further improvement of a minimax theorem of
Borenshtein and Shul'man, {\it J. Nonlinear Convex Anal.}, 2:279-283, 2001.\par
\smallskip
\noindent
[5]\hskip 5pt B. Ricceri, Minimax theorems for limits of 
parametrized functions
having at most one local minimum lying in a certain set, {\it Topology Appl.}, 
153:3308-3312, 2006.\par
\smallskip
\noindent
[6]\hskip 5pt B. Ricceri, Recent advances in minimax theory and
applications, in ``Pareto Optimality, Game Theory and Equilibria'', 
A. Chinchuluun, P.M. Pardalos, A. Migdalas, L. Pitsoulis eds., Springer, 23-52,
2008.\par
\smallskip
\noindent
[7]\hskip 5pt B. Ricceri, 
{\it Nonlinear eigenvalue problems},  
in ``Handbook of Nonconvex Analysis and Applications'' 
D. Y. Gao and D. Motreanu eds., International Press, 543-595, 2010.\par
\smallskip
\noindent
[8]\hskip 5pt B. Ricceri, A strict minimax inequality criterion and some of its consequences, {\it Positivity}, 16:455-470,
2012.\par
\smallskip
\noindent
[9]\hskip 5pt B. Ricceri, Energy functionals of Kirchhoff-type problems having multiple global minima, {\it Nonlinear Anal.}, 115:130-136,
2015.\par
\smallskip
\noindent
[10]\hskip 5pt B. Ricceri, A minimax theorem in infinite-dimensional topological vector spaces, preprint.\par
\smallskip
\noindent
[11]\hskip 5pt J. Saint Raymond, Connexit\'e des sous-niveaux des fonctionnelles int\'egrales, {\it Rend. Circ. Mat. Palermo}, 
44:162-168, 1995.\par
\smallskip
\noindent
[12]\hskip 5pt  E. Zeidler,  Nonlinear functional analysis and
its applications, vol. III, Springer-Verlag, 1985.\par

\bye

\bye